\documentclass[reqno]{amsart}

\usepackage{amsmath,amsthm,amssymb,amscd,enumerate,mathtools,enumitem}
\usepackage{amsfonts}
\usepackage{rotating}
\usepackage{euscript}
\usepackage{pst-node}
\usepackage{epsfig,verbatim}
\usepackage{stackengine,scalerel}
\stackMath
\newcommand\qfrac[3][1pt]{\frac{%
  \ThisStyle{\addstackgap[#1]{\SavedStyle#2}}}{%
  \ThisStyle{\addstackgap[#1]{\SavedStyle#3}}%
}}

\def\be{\begin{equation}}
\def\ee{\end{equation}}

\def\C{{\mathbb C}} 
\def\f{\EuScript}
 
\def\P{{\mathbb P}}

\def\T{{\mathbb T}}

\def\R{{\mathbb R}}

\def\phi{{\varphi}}

\def\deg{{\rm deg\,}}

\def\Sing{\operatorname{Sing}}
\def\Im{\operatorname{Im}}
\def\bp{\begin{proposition}}
\def\ep{\end{proposition}}

\def\bt{\begin{theorem}}
\def\et{\end{theorem}}
\def\br{\begin{remark}}
\def\er{\end{remark}}
\def\be{\begin{equation}}
\def\bee{\begin{equation*}}
\def\l{\label}

\def\ee{\end{equation}}
\def\eee{\end{equation*}}
\def\bl{\begin{lemma}}
\def\el{\end{lemma}}
\def\bc{\begin{corollary}}
\def\ec{\end{corollary}}
\def\pr{\noindent{\it Proof. }}

\def\bd{\begin{definition}}
\def\ed{\end{definition}}

\def\cC{{\mathcal C}}

\def\fR{\EuScript{R}}
\def\fM{\EuScript{M}}

\def\bC{\widehat{\C}} 
\def\bR{\widehat{\R}} 
\def\E{\widehat{L}}   

\newtheorem{theorem}{Theorem}[section]
\newtheorem{lemma}[theorem]{Lemma}
\newtheorem{definition}[theorem]{Definition}
\newtheorem{corollary}[theorem]{Corollary}
\newtheorem{proposition}[theorem]{Proposition}
\newtheorem{problem}[theorem]{Problem}

\theoremstyle{definition}

\theoremstyle{definition}
\newtheorem{remark}[theorem]{Remark}
\def\bpr{\begin{problem}}
\def\epr{\end{problem}}

\mathtoolsset{showonlyrefs}

\begin{document}

\title[On intersection of lemniscates   of rational functions]
      {On intersection of lemniscates\\ of rational functions}

\author[S. Orevkov]{Stepan Orevkov}

\address{IMT, l'universit\'e Paul Sabatier, 119 route de Narbonne,
         31062 Toulouse, France}
\address{Steklov Mathematical Institute, ul. Gubkina 8, Moscow 119991, Russia}

\email{orevkov@math.ups-tlse.fr}

\author[F. Pakovich]{Fedor Pakovich}
\address{Department of Mathematics, 
Ben Gurion University of the Negev, P.O.B. 653, Beer Sheva,  8410501, Israel}
\email{pakovich@math.bgu.ac.il}

\begin{abstract}
For a non-constant complex rational function $P$, the  {\it lemniscate}  of $P$ is
defined as the set of points $z\in \C$ such that $\vert P(z)\vert =1$.
The lemniscate of $P$ coincides with the set of real points of the algebraic
curve given by the equation $L_P(x,y)=0$,  where  $L_P(x,y)$ is 
the numerator of the rational function $P(x+iy)\overline{ P}(x-iy)-1.$  
In this paper, we study the following two questions: under what conditions two
lemniscates have a common component, and under what conditions the algebraic curve
$L_P(x,y)=0$ is irreducible. 
In particular, we provide a sharp bound for the number of complex solutions of
the system $\vert P_1(z)\vert =\vert P_2(z)\vert =1$,
where $P_1$ and $P_2$ are rational functions.   
\end{abstract}


\maketitle

\section{Introduction}
Let $P$ be a non-constant complex rational function on the Riemann sphere $\bC$.
The \emph{lemniscate} of $P$ is defined as
\be
                 \f L_P=\{z\in \C: \, \vert P(z)\vert =1\}.
\ee
The paper is devoted to two more or less independent problems concerning the lemniscates.
The first problem is to determine the maximal number of intersection points of $\f L_{P_1}$ and $\f L_{P_2}$,
i.e., the number of solutions of the system of equations
\be
   \label{sys1}
                 \vert P_1(z)\vert =\vert P_2(z)\vert =1,
\ee
where $P_1$ and $P_2$ are rational functions of given degrees $n_1$ and $n_2$.
The second problem is to find out when a lemniscate (considered as a real algebraic curve in $\R^2$)
is irreducible.

Our result about the intersection of lemniscates is the following statement.

\begin{theorem}\l{3}
Let $P_1$ and $P_2$ be non-constant complex rational functions
of degrees $n_1$ and $n_2$.
The following three conditions are equivalent:
\begin{enumerate}
\item[{\rm(i)}]
       ${\f L}_{P_1}$ and ${\f L}_{P_2}$ have more
       than $2n_1 n_2$ common points.
\item[{\rm(ii)}]
       ${\f L}_{P_1}\cap{\f L}_{P_2}$ is infinite.
\item[{\rm(iii)}]
       $P_1=B_1\circ W$ and $P_2=B_2\circ W$ for some rational functions
       $W$, $B_1$, $B_2$ such that
       each of $B_1,B_2$ maps the unit circle $\T$ to itself.
\end{enumerate}
 Furthermore, for any natural $n_1$ and $n_2$ 
there exist rational functions of degrees $n_1$ and $n_2$ such that
$\f L_{P_1}$ and $\f L_{P_2}$ 
have exactly $2n_1n_2$ intersection points.  
\end{theorem}

We remark that
for a rational function $B$, the condition $B(\T)\subset\T$
can be written in the form
\be
    \l{eq.bla}
               B(z)\overline{B}(1/z)=1 \qquad\text{for all $z\in\bC$.}
\ee
Such functions are known in the complex analysis under the name of
{\it quotients of finite Blaschke products}.

In a recent paper \cite{ps}, the equivalence of conditions $\rm(i), \rm(ii), \rm(iii)$
was proven with a weaker bound $(n_1+n_2)^2$ instead of the bound $2n_1n_2$.
Notice that the result of \cite{ps} implies
the main result of the paper \cite{AR} by Ailon and Rudnick, which  is equivalent to
the following statement: if $P_1$ and $P_2$ are complex polynomials, then 
\be
     \#\bigcup_{k=1}^\infty \{z \in \C~:~P_1(z)^k = P_2(z)^k =1\} \le C(P_1,P_2),
\ee
for some   constant $C(P_1,P_2)$ that depends only on $P_1$ and $P_2$,
unless for some non-zero integers $m_1$ and $m_2$ the equality 
\be
    \l{PP} 
                   P_1^{m_1}(z)P_2^{m_2}(z)=1
\ee
holds (see Section \ref{sect.ipl} below).
In addition, the result of  \cite{ps} answers the question of 
Corvaja, Masser, and Zannier (\cite{CMZ}) about the intersection of an irreducible
curve $\cC$ in $\C^*\times \C^*$ with $\T\times \T$,  in  case $\cC$ has genus zero
and is parametrized by rational functions $P_1,$ $P_2$  
(see \cite{ps} for more details). These applications of \cite{ps} stem from
the fact that  the numbers  $C(P_1,P_2)$ and $\#(\cC\cap \T\times \T)$ 
obviously are bounded from above by the number of solutions of \eqref{sys1}.
Thus, a sharp bound for the last number is of great interest,
and our Theorem~\ref{3} provides it. 

In brief, our proof of Theorem~\ref{3} given in Section 2 goes as follows.
If $P$ is a complex rational function of degree $n$, then under the standard
identification of $\C$ with $\R^2$ the lemniscate $\f L_P$
coincides with the set of real points of the affine algebraic curve of degree $2n$
given by the equation $L_P(x,y)=0$,  where $L_P(x,y)$ is 
the numerator of the rational function
\be
   \l{eq.L}
                P(x+iy)\overline{ P}(x-iy)-1.
\ee
After the linear change of variables  $z=x+iy$, $w=x-iy$ (in $\C^2$) the Newton
polygon of $L_P$ becomes the square $n\times n$. Thus, the B\'ezout Theorem for
bihomogeneous polynomials (\cite[\S 4.2.1]{shaf}), which is also a simplest case of the
Bernstein-Kushnirenko Theorem \cite{b}, implies that if (i) holds,
then $L_{P_1}$ and $L_{P_2}$ have a common factor, i.e.,
the system $L_{P_1}=L_{P_2}=0$ has infinitely many {\it complex} solutions.
This is not yet (ii), but we prove a kind of ``real'' version of the 
B\'ezout theorem (Proposition \ref{prop.bezout}),
which implies in particular that if (i) holds, then
the system $L_{P_1}=L_{P_2}=0$ has infinitely many {\it real} solutions
(Corollary \ref{cor.bezout}(c)).
This gives us the implication
(i) $\Rightarrow$ (ii). 
In turn, the implication (ii) $\Rightarrow$ (iii) is deduced from
a ``rational'' version of 
the Cartwrite Theorem \cite{c} (Corollary~\ref{car2}). 
Finally,  since (iii) implies that $\f L_W\subset\f L_{P_1}\cap\f L_{P_2}$,
the   implication (iii) $\Rightarrow$ (i) is obvious.

Notice that the equivalence of conditions $\rm(i), \rm(ii), \rm(iii)$ can be proven also by modifying the approach of the paper \cite{ps} (see Section \ref{23}).

Another problem considered in this paper (in Section~\ref{s2}) is the following.
Given  a rational function $P$, under what conditions the curve $L_{P}(x,y)=0$ is irreducible over $\C$?  
It is not hard to see that the following ``composition condition'' is sufficient
for the reducibility: there exists a quotient of finite Blaschke products $B$
of degree at least 2 and a rational function $W$ such that $P=B\circ W$.
Notice that if $P$ is a {\it polynomial}, the composition condition
reduces to the condition that $P=W^k$ for some polynomial $W$ and $k\ge 2$.
Its necessity  for the reducibility of $L_P(x,y)$ in the polynomial case 
was established by the first author in \cite{o} (notice that this result has
found applications to complex dynamics, see \cite{ok}).
In Section~\ref{sect.ce},
we show however that for {\it rational} $P$ the reducibility of
$L_P(x,y)=0$ does not imply in general the composition condition.

Our approach to the problem of irreducibility of algebraic curves
$L_P(x,y)=0$ is based on the observation that a change of variables allows us to consider this problem
 in the context of the more general problem of irreducibility
of ``separated variables" curves $P(x)-Q(y)=0$, where $P(z)$ and $Q(z)$ are rational functions. In particular, we deduce our examples of reducible $L_P(x,y)$ for $P$ not satisfying the composition condition from examples of reducible separated variables curves found in \cite{cc}.
On the other hand, modifying the arguments from \cite{o}, we provide a handy
{\it sufficient} condition for the irreducibility of separated variables curves in case  
{\it one} of rational functions $P(z)$ and $Q(z)$ is a polynomial (Theorem~\ref{tp}).


\section{Intersection of lemniscates}\label{sect.intersect}

\subsection{The Cartwright theorem}\label{sect.cart}

 For a simple closed curve $\Gamma$ in $\C$, we denote by $\fM_{\Gamma}$ the set of all non-constant functions 
meromorphic on $\C$ having modulus one on $\Gamma$.
The following result was proved by Cartwright in \cite{c}. 

\begin{theorem}\l{car}
 Let  $\Gamma$ be a simple closed curve in $\C$ such that $\fM_{\Gamma}$ is non\-empty. Then there exists a  function $\varphi\in \fM_{\Gamma}$ such that 
each function in $\fM_{\Gamma}$ may be written in the form $f=B\circ \varphi$, where $B$ is a quotient of finite Blaschke products.
\end{theorem}

A detailed discussion and generalizations of the  Cartwright theorem can be found in the papers  \cite{stsu}, \cite{ste}. 
Below we need the following specialization of the  Cartwright theorem, where    
the notation $\fR_{\Gamma}$ stands for the set of all non-constant complex rational functions on $\C$ having modulus one on $\Gamma$.

\begin{corollary}\l{car2}
Let  $\Gamma$ be a simple closed curve in $\C$ such that $\fR_{\Gamma}$ is non\-empty.
Then there exists a rational function $W\in \fR_{\Gamma}$ such that 
each function in $\fR_{\Gamma}$ may be written in the form $P=B\circ W$,
where $B$ is a quotient of finite Blaschke products.
\end{corollary}

\pr Applying Theorem \ref{car}, we conclude that there exists a meromorphic function $W\in \fM_{\Gamma}$ such that if $P\in \fR_{\Gamma}$, then $P=B\circ W$ for some  quotient of finite Blaschke products $B$. However, 
 since the great Picard theorem implies that any non-rational   function $W$ meromorphic  on $\C$ takes all but at most two values in $\mathbb C\mathbb P^1$ infinitely often, this equality implies  that $W$ is rational. \qed


\subsection{A real version of the B\'ezout Theorem.}\l{sect.bezout}
Let us recall that the classical B\'ezout theorem about intersections of curves in
$\C\P^2$ implies that the number of intersection points of two affine algebraic curves
$F(x,y)=0$ and $G(x,y)=0$ of degrees $m$ and $n$ in $\C\times \C$ does not exceed $mn$,
unless the polynomials $F(x,y)$ and $G(x,y)$ have a non-constant common factor in
$\C[x,y]$.
  In case the bidegrees $(m_1,m_2)$ and $(n_1,n_2)$ of $F(x,y)$ and $G(x,y)$ are relatively small with respect to $m$ and $n$, a better bound can be obtained from the 
bihomogeneous B\'ezout theorem, 
which implies that  the number of intersection 
points of $F(x,y)=0$ and $G(x,y)=0$   does not exceed $n_1m_2+n_2m_1$, unless
$F(x,y)$ and $G(x,y)$ have a non-constant common factor (see \cite[\S4.2.1]{shaf}).
In the proof of Theorem \ref{3}, we will use a real version of this last bound provided
by Corollary~\ref{cor.bezout}(c) below.

Let $X$ be a compact non-singular complex algebraic (or analytic) variety.
A {\it real structure} on $X$ is an anti-holomorphic involution
$\sigma:X\to X$. We denote the variety $X$ endowed with the real structure $\sigma$
by $X^\sigma$. A point $p$ of $X^\sigma$ is called {\it real} if $\sigma(p)=p$.
The set of real points of $X^\sigma$ (the {\it real locus} of $X^\sigma$) is
denoted by $\R X^\sigma$. A basic example is a projective variety in $\C\P^n$
defined by polynomial equations with real coefficients endowed with the
involution of complex conjugation. In this case $\R X^\sigma$ is the subset of
$\R\P^n$ defined by the same equations.

If $X=\C\P^1\times\C\P^1$, then
there are two non-isomorphic real structures on $X$:
$$
   \sigma_h:(x,y)\mapsto(\bar x,\bar y), \qquad
   \sigma_e:(x,y)\mapsto(\bar y,\bar x).
$$
The real loci $\R X^{\sigma_h}$ and $\R X^{\sigma_e}$ are, respectively,
$\R\P^1\times\R\P^1$ (a torus) and the image of $\C\P^1$ in $X$ under
the embedding $z\mapsto(z,\bar z)$ (a sphere).
Note that $X^{\sigma_h}$ and $X^{\sigma_e}$ are isomorphic
to the complexifications of, respectively, hyperboloid and ellipsoid in $\C\P^3$
endowed with the usual complex conjugation:
$$
   X^{\sigma_h}\cong\{z\mid z_0^2+z_1^2=z_2^2+z_3^2\}, \qquad
   X^{\sigma_e}\cong\{z\mid z_0^2=z_1^2+z_2^2+z_3^2\};
$$
here $z$ stands for $(z_0:z_1:z_2:z_3)$.
Indeed, it is straightforward to check that the following mapping
$X\to\C\P^3$ defines the required isomorphism:
$$
  (x,y)\mapsto
    (x_0y_0+x_1y_1:\alpha(x_0y_1-x_1y_0):x_0y_0-x_1y_1:x_0y_1+x_1y_0),
$$
where $\alpha=1$ for $\sigma_h$ and $\alpha=i$ for $\sigma_e$;
here $x$ and $y$ stand for $(x_0:x_1)$ and $(y_0:y_1)$.

Let us recall that
if $C$ is an irreducible curve on a smooth compact algebraic surface $X$, then
by the genus formula (see e.g.~\cite[\S4.4.1]{shaf}) the number of singular
points of $C$ is bounded from above by its {\it arithmetic genus} $p_a(C)$,
which in this case can be computed by the formula
$$
     2p_a(C) = 2 + C\cdot(C+K_X),
$$
where $K_X$ is the canonical class of $X$.

\begin{proposition}\l{prop.bezout}
Let $X$ be a smooth compact complex algebraic surface
endowed with a real structure $\sigma$.
Let $A$ and $B$ be $\sigma$-invariant algebraic curves on $X$.
Suppose that $\R A\cap\R B$ is finite and $C^2\ge p_a(C)$
for each irreducible component $C$ of $A\cup B$. Then
\be
      \#(\R A\cap\R B) \le A\cdot B.
\ee
\end{proposition}

\pr
We say that a $\sigma$-invariant curve $C$ on $X$ is {\it $\sigma$-irreducible} if
$C=D\cup\sigma(D)$ where $D$ is irreducible.

Let us first consider the case when $A$ and $B$ are $\sigma$-irreducible.
If $A\ne B$, then $A$ and $B$ have no common components, and hence
$$
    \#(\R A\cap\R B) \le \#(A\cap B) \le A\cdot B.
$$
Suppose now that $A=B$.
Let us show that $\R A\subset\Sing(A)$, where $\Sing(A)$
is the set of singular points of $A$.
Indeed, for any $p\in\R A$ one can choose local coordinates $(z,w)$ in
a neighborhood of $p$ such that $\R X=\{\Im z=\Im w=0\}$ and
$A=\{f(z,w)=0\}$, where $f$ is a polynomial with real coefficients.
Thus, if $A$ were non-singular at $p$, then $\R A$ would be infinite by
the Implicit Function Theorem,
which contradicts our hypothesis that $\R A\cap\R B$ is finite.
Thus, $\R A\subset\Sing(A)$.

It follows that if $A=B$ is irreducible (in the usual sense), we have
$$
    \#(\R A\cap\R B)=\#\R A \le \#\Sing(A) \le p_a(A) \le A^2 = A\cdot B.
$$
On the other hand,
if $A=B$ is reducible, then $A=C\cup\sigma(C)$ where $C$ is
irreducible and $\sigma(C)\ne C$.
Therefore, we have $\R A\subset C\cap\sigma(C)$, and hence
$$
  \#(\R A\cap\R B) = \#\R A \le C\cdot\sigma(C)
      \le C^2+\sigma(C)^2+2C\cdot\sigma(C) = A^2 = A\cdot B. 
$$
This completes the proof in the case where $A$ and $B$ are $\sigma$-irreducible.

Now we proceed to the general case.
Let $A=A_1\cup\dots\cup A_k$ and $B=B_1\cup\dots\cup B_l$
where each $A_i$ and each $B_j$ is $\sigma$-invariant and $\sigma$-irreducible.
Then
$$
   \#(\R A\cap\R B) \le \sum_{i,j}\#(\R A_i\cap\R B_j)
                     \le \sum_{i,j} A_i\cdot B_j = A\cdot B. \eqno{\Box}
$$

\medskip

\bc\l{cor.bezout}
{\rm(a).} 
Let $F(x_1,x_2)$ and $G(x_1,x_2)$ be polynomials with real coefficients of degrees
$m$ and $n$.
Then the number of real solutions of the system $F=G=0$ is either infinite
or bounded above by $mn$.

\smallskip{\rm(b).}
Let $F(x_1,x_2)$ and $G(x_1,x_2)$ be polynomials with real coefficients such that
$\deg_{x_k} F=m_k$ and $\deg_{x_k} G
=n_k$, $k=1,2$.
Then the number of real solutions of the system $F=G=0$ is either infinite
or bounded above by $m_1n_2+m_2n_1$.

\smallskip{\rm(c).}
Let $F(z,w)$ and $G(z,w)$ be polynomials with complex coefficients such that
$F(z,w) = \overline{F}(w,z)$ and $G(z,w) = \overline{G}(w,z)$. Let
$\deg_z F=m$ and $\deg_z G=n$.
Then the number of solutions of the system $F=G=0$ belonging to the set
\hbox{$\{(z,w)\mid w=\bar z\}$}  is either infinite or bounded above by $2mn$.
\ec

\pr
(a). We apply Proposition~\ref{prop.bezout}
to the curves in $\C\P^2$ (with the standard real structure) defined by
the homogeneous equations
$$
   x_0^m F(x_1/x_0,x_2/x_0) = 0, \qquad x_0^n G(x_1/x_0,x_2/x_0) = 0,
$$
and observe that, for a curve $C$ of degree $d$ in $\C\P^2$, we have
(see \cite[\S4.2.3]{shaf})
$$
    p_a(C) = (d-1)(d-2)/2 \le d^2 = C^2.
$$
(b). We apply Proposition~\ref{prop.bezout} to $X=\C\P^1\times\C\P^1$
endowed with the real structure $\sigma_h$ (see above) and to the curves
defined by the bihomogeneous equations
$$
   u_1^{m_1}u_2^{m_2} F(x_1/u_1,x_2/u_2) = 0, \qquad
   u_1^{n_1}u_2^{n_2} G(x_1/u_1,x_2/u_2) = 0.
$$
In this case, for a curve $C$ of bidegree $(d_1,d_2)$ in $X$, we have
(see \cite[\S4.2.3]{shaf})
$$
   p_a(C) = (d_1-1)(d_2-1) \le 2d_1d_2 = C^2.
$$
(c). The same proof as in (b), but with $\sigma_e$ instead of $\sigma_h$.
(Note that our hypothesis about $F$ and $G$ implies that their bidegrees
are $(m,m)$ and $(n,n)$ respectively.)
\qed

\medskip
Notice that 
higher-dimensional analogs of Corollary~\ref{cor.bezout} do not hold.
Indeed,
let $F_1=P(x)+P(y)$, $F_2=F_3=z$, where $P(x)=\prod_{k=1}^n(x-k)^2$, $n\ge 3$
(see \cite[Example 13.6]{ful}).
Then the number of real solutions of the
system of equations $F_1=F_2=F_3=0$ is finite but greater than $\prod\deg F_i$.


\subsection{Proof of Theorem~\ref{3}}\label{23}
The implication (i) $\Rightarrow$ (ii) immediately
follows from Corollary~\ref{cor.bezout}(c),
because the embedding $\C\to\C^2$, $z\mapsto(z,\bar z)$ identifies
$\f L_{P_k}$ with
$$
    \{(z,w)\,:\,P_k(z)\overline P_k(w)=1\}\,\cap\, \{(z,w)\,:\,w=\bar z\},
    \quad k=1,2.
$$

\smallskip
The implication (ii) $\Rightarrow$ (iii) follows from Corollary~\ref{car2}. Indeed,
suppose that the intersection $\f L_{P_1}\cap\f L_{P_2}$ is infinite.
Then the complexifications of $\f L_{P_1}$ and $\f L_{P_2}$ have a common real
component $A$ whose real locus $\R A$ is also infinite.
By composing $P_1$ and $P_2$ with a M\"obius transformation,
we may assume that $|P_j(\infty)|\ne 1$, \hbox{$j=1,2$}, so that $\R A$ is compact.
Recall that the real locus of a real algebraic curve in neighborhood of every its 
non-isolated point is homeomorphic to  a ``star'' with an even number of branches
(see e.g.~\cite{gh}, p.~104). Thus, the set of non-isolated points of $\R A$
is homeomorphic to a graph in $\R^2$, all of whose vertices have even valency. Such a graph necessarily have a cycle $D$, and, by construction, the rational functions $P_1$ and $P_2$ have modulus one on $D$.
Applying now Corollary~\ref{car2} to $D$, we obtain (iii).

\smallskip
To prove (iii) $\Rightarrow$ (i), it is enough to observe that if $B$
 is a quotient of finite Blaschke products, then 
$B^{-1}(\T)$ contains $\T$. Thus, (iii) implies that both
$\f L_{P_1}$ and $\f L_{P_2}$ contain the infinite set 
 $\f L_W$ as a subset.  

\smallskip
To prove the last part of the theorem, let us fix a M\"obius transformation $\nu$ that
maps the real line to the unit circle, and observe that the lemniscate of
$P_1=\nu\circ z^{n_1}$ is a union of $n_1$ lines on $\C$ passing through the origin.
Under the identification of $\bC$ with $S^2$ via the stereographic projection,
$\f L_{P_1}$ becomes a union of $n_1$ big circles passing through two
antipodal points $a_1,b_1$. Let $\delta$ be a M\"obius transformation of $\bC$
corresponding to an isometry of $S^2$, and $P_2=\nu\circ z^{n_2}\circ \delta$.
Then $\f L_{P_2}$ is a union of $n_2$ big circles passing through two antipodal points
$a_2,b_2$. Any two distinct big circles intersect at two points. Hence, if $\delta$
is chosen generically (so that
$\{a_1,b_1\}\cap\f L_{P_2}=\{a_2,b_2\}\cap\f L_{P_1}=\varnothing$), then
$\#(\f L_{P_1}\cap\f L_{P_2})=2n_1 n_2$.
\qed 

\medskip

Note that Theorem~\ref{3} can be proved somewhat shorter using the approach of the paper \cite{ps},  
where instead of the intersection $\f L_{P_1}\cap\f L_{P_2}$ an algebraic curve $\cC$ parametrized by  $P_1$ and $P_2$ is considered. Under this approach, the proof of Theorem~\ref{3} reduces to 
proving that 
if an algebraic curve $\cC$  parametrized by rational functions $P_1$ and $P_2$ of degrees $n_1$ and $n_2$ has more than $2n_1n_2$ points whose coordinates have modulus one, 
then it can be parametrized by some quotients of finite Blaschke products. On the other hand, considering instead of the functions 
$P_1$, $P_2$ the functions 
$$\widehat P_1=T\circ P_1\circ T^{-1}, \ \ \ \ \widehat P_2=T\circ P_2\circ T^{-1},$$ 
where 
$$T(z) = i\frac{1+z}{1-z},$$ 
one can easily see that the last statement is equivalent to the statement that    
 if an algebraic curve $\cC$  paramet\-ri\-zed by rational functions $P_1$ and $P_2$ of degrees $n_1$ and $n_2$  has more than $2n_1n_2$ real points, 
then it can be paramet\-ri\-zed by some rational functions with real coefficients. Finally, the last statement follows from 
the bihomogeneous B\'ezout Theorem since if  $\cC$ is not defined over $\R$, then the set 
$$\cC\cap \R^2=\cC\cap \overline \cC$$ is finite and  contains at most $\cC\cdot \overline \cC=2n_1n_2$ points (in \cite{ps} the   usual version of the B\'ezout Theorem was used).

Our proof of Theorem~\ref{3}, while longer, considers the lemniscates directly  and relates the initial problem with  
a real version of the B\'ezout Theorem, which might have an independent interest.


\subsection{Intersection of polynomial lemniscates}\l{sect.ipl}

The polynomial version of Theorem \ref{3} is the following statement.

\bt \l{pol}
Let $P_1$ and $P_2$ be non-constant complex polynomials of degrees $n_1$ and $n_2$.
The following three conditions are equivalent:
\begin{enumerate}
\item[{\rm(i)}\;\;]
       ${\f L}_{P_1}$ and ${\f L}_{P_2}$ have more
       than $2n_1 n_2$ common points.
\item[{\rm(ii)}\;\;]
       ${\f L}_{P_1}\cap{\f L}_{P_2}$ is infinite.
\item[{\rm(iii)*}]
    $P_1=P^{n_1},$ $P_2=c P^{n_2}$
for some polynomial $P$, natural $n_1,n_2$,
and $c\in \C$ with $\vert c\vert =1$.
\end{enumerate}
\et
 
\pr It is clear that Condition (iii)* is a particular case of Condition (iii)
for $B_1=z^{n_1},$ $B_2=cz^{n_2},$ and $W=P.$ Thus, in view of Theorem~\ref{3}, it is enough to show that if $P_1$ and $P_2$ are polynomials, then (iii) reduces to (iii)*.

To prove the last statement, let us observe that if $P_1$ and $P_2$ are polynomials,
then each of the functions $B_1$ and $B_2$ has a unique pole,
and this pole is common for $B_1$ and $B_2$.
By condition \eqref{eq.bla}, each pole of $B_k$, $k=1,2$,
is symmetric with respect to $\T$ to a zero of the same multiplicity and vice versa.
Therefore, there exists $a\in \bC\setminus \T$ such that
$$
           B_1=c_1\Big( \frac{z-a}{1-\overline{a}z}\Big)^{n_1}, \qquad
           B_2=c_2\Big( \frac{z-a}{1-\overline{a}z}\Big)^{n_2}
$$
for some $n_1,n_2\geq 1$ and $c_1,c_2\in \C$ with $\vert c_1 \vert =\vert c_2 \vert=1$, implying that  
\be
   \l{b1}      
    P_1=c_1 \widetilde P^{n_1}, \qquad
    P_2=c_2 \widetilde P^{n_2},
\ee
where
$$
   \widetilde P=\frac{z-a}{1-\overline{a}z}\circ W.
$$
Finally, it is easy to see that $\widetilde P$ is a polynomial and (iii)* holds for  
$P=c_1^{\frac{1}{n_1}}\widetilde P$ and $c=c_2c_1^{-{n_2}/{n_1}}$.
\qed
\medskip

Notice that in the paper \cite{ps} it is erroneously  claimed that in the polynomial
case Condition (iii) of Theorem~\ref{3} simply reduces to the condition that
$P_1 = P^{m_1}$, $P_2 = P^{m_2}$ for some polynomial $P$.
This inaccuracy however does not affect the main application of results
of \cite{ps} in the polynomial case:
the result of Ailon and Rudnick mentioned in the introduction.
Indeed, if $c$ in (iii)* is not a root of unity, then the system
\be
   \l{sim}
                 P_1(z)^k =1, \ \ \  P_2(z)^k =1
\ee
has no solutions for any $k\geq 1$, since
(iii)* and \eqref{sim} imply that
$$
     P^{n_1n_2k}(z)=c^{n_1k}P^{n_1n_2k}(z)=1.
$$
On the other hand, if $c^l=1$, then \eqref{PP} holds for
$m_1=n_2l$ and $m_2=-n_1l$.


\subsection{On sharpness of the bound $2n_1n_2$ in the polynomial case}\l{sect.sharp}

The last statement of Theorem~\ref{3} states that the bound $2n_1n_2$ is sharp
if we speak of all rational functions. However, it does not seem so when we restrict
ourselves to polynomials only. The maximal number of intersection points of
two polynomial lemniscates that we succeeded to realize, is given in
the following statement.

\bp\l{exa.poly}
Let $1\le n_1\le n_2$ and $d=\gcd(n_1,n_2)$.
Then there exist polynomials $P_1$ and $P_2$
such that $\deg P_k=n_k$, $k=1,2$, and
$$
    \#(\f L_{P_1}\cap\f L_{P_2}) = M(n_1,n_2)
    := n_1 n_2 + n_2 + de,\qquad e
    = \begin{cases} 0, &\text{if $n_1/d$ is odd,}\\
                    1, &\text{if $n_1/d$ is even.}
      \end{cases} 
$$
\ep

\pr
Let $P_k(z) = (z/r_k)^{n_k} - 1$, $k=1,2$, where $r_k>0$.
The lemniscate $\f L_{P_k}$ is ``flower-shaped'', i.e. it is a union of
$n_k$ loops (``petals'') outcoming from $0$.
It is clear that $\f L_{P_k}$ tends to a union of $n_k$ lines
(we denote it by $\f L_0$) when $r_k\to\infty$.
Let us fix $r_2$.
It is not difficult to show that for a suitably chosen rotation $R$ we have
$$
    \#\big(\f L_0\cap R(\f L_{P_2})\setminus\{0\}) = n_2+de.
$$
A further small shift of $R(\f L_{P_2})$
produces $n_1 n_2$ additional crossings near $0$.
Finally, we approximate $\f L_0$ by $\f L_{P_1}$ by choosing $r_1\gg r_2$.
\qed

\medskip
Note that the number $M(n_1,n_2)$ in Proposition~\ref{exa.poly}
coincides with the upper bound $2n_1 n_2$ given by Theorem~\ref{3}
if and only if $n_1=1$.


\section{Irreducibility of lemniscates}\l{s2} 
\subsection{Irreducibility of separated variables curves and lemniscates} 
Let us recall that a {\it separated variable curve} is an algebraic curve in $\C^2$
given by the equation $E_{P,Q}(x,y) = 0$, where $E_{P,Q}(x,y)$ is the numerator of
$P(x) - Q(y)$ for some non-constant complex rational functions $P(z)$ and $Q(z)$.
 The irreducibility problem for separated variable curves is quite old and still
 widely open (see \cite{fc} for an introduction to the subject). One of the most
 complete results in this direction is a full classification of reducible curves
 $E_{P,Q}(x,y)=0$ in the case when $P$ and $Q$ are indecomposable polynomials. 
 All possible ramifications of such $P$ and $Q$ were described by Fried in \cite{f4},
 and polynomials themselves were found by Cassou-Nogu\`es and Couveignes in \cite{cc}.
For further progress, we refer the reader to the recent paper \cite{nef2} and the
bibliography therein. Here and below, by the irreducibility we always mean the
irreducibility over $\C$. 

Let us recall that $L_P(x,y)$ is the numerator of the rational function \eqref{eq.L},
and set
$$
      \E_P(x,y) = E_{P,\,1/\,\overline P}(x,y).
$$
The irreducibility problem for curves $L_P(x,y)=0$ defining lemniscates is a particular case of the irreducibility problem for separated variables curves due to the following
statement, which is immediate from the identity
$L_P(x,y) = \E_P(x+iy,x-iy)$.

\bl \l{u1}
 Let $P$ be a non-constant complex rational function. Then the curve
$L_P(x,y)=0$ is irreducible if and only if the curve
$\E_P(x,y)=0$ is irreducible. \qed 
\el

 The following theorem was proved   in the paper \cite{o}.

\bt  \l{or} Let $P$ and $Q$ be non-constant complex polynomials.
 Then the curve  $E_{P,\,1/Q}(z,w)=0$ is reducible
if and only if
$$
    P(z)=P_1(z)^d \qquad\text{and}\qquad Q(w)=Q_1(w)^d
$$
for some $d>1$ and polynomials $P_1(z)$ and $Q_1(w)$.
\et

It is easy to see that combined with Lemma \ref{u1},  Theorem \ref{or} implies the following result also proved in \cite{o}. 

\bt \l{or2}
Let $P$ be a non-constant  complex polynomial. Then the polynomial $L_P(x,y)$ is reducible if and only if $$
    P(z)=P_1(z)^d 
$$
for  some $d>1$ and polynomial $P_1(z)$. 
\et

Notice that
since $|P_1(z)^d|=1\Leftrightarrow |P_1(z)|=1$, this theorem implies that
any polynomial lemniscate is the zero set of an irreducible polynomial in $x$ and $y$.

\bd\l{def.cc}
{\rm A complex rational function $P(z)$ satisfies the
\emph{Composition Condition} if
there exist a quotient of finite Blaschke products $B$ of degree at least two
and a non-constant rational function $W$ such that $P=B\circ W$. }
\ed

Theorem \ref{or2} shows that the Composition Condition is necessary and sufficient
for the reducibility of $L_P(x,y)$ in the polynomial case.
Moreover, it is easy to see that the Composition Condition is {\it sufficient}
for the reducibility of $L_P(x,y)$ for any rational function $P$. Indeed,
it follows from \eqref{eq.bla} that for any quotient of finite Blaschke products $B$
and a rational function $W$ the equality $W(x)=1/\,\overline{W}(y)$ 
for some $x,y\in \C$ implies the equality  
$$
    P(x)=B(W(x))=\qfrac[1pt]{1}{\overline{B}(1/W(x))}
        =\qfrac[1pt]{1}{\overline{B}(\overline{W}(y))}
        =\qfrac[1pt]{1}{\overline{P}(y)}.
$$
Therefore, the curve $\E_W(x,y)=0$ is a component of the curve $\E_P(x,y)=0$,
implying that the latter curve is reducible since the considered curves 
have different degrees (in view of the assumption $\deg B>1$). Thus, the curve $L_P(x,y)=0$ is reducible by Lemma~\ref{u1}.

In the next section, we show that Composition Condition is not necessary
for the reducibility of $L_P$, while in this section, 
modifying the idea of \cite{o}, we establish a sufficient condition 
for the irreducibility of an algebraic curve $E_{P,Q}(x,y)=0$ in the case when
one of the rational functions $P$, $Q$ is a polynomial.

Namely, we prove the following result. 


\bt \l{tp} Let $P$ be a polynomial of degree $n\geq 1$, and $Q$ a rational function.
Assume that multiplicities $q_1,q_2,\dots,q_l$ of poles of $Q$ satisfy the
condition  $\gcd(q_1,q_2,\dots,q_l,n)=1$. 
Then the curve $E_{P,Q}(x,y)=0$ is irreducible.  
\et

Notice that like Theorem \ref{or}, Theorem \ref{tp} easily implies Theorem \ref{or2}.

\medskip
\noindent {\it First proof of Theorem \ref{tp}}
(cf.~the proof of \cite[Theorem 1]{o}).
Let $C$ be the closure of the curve $E_{P,Q}=0$
in $\C\P^1\times\C\P^1$, that is, the curve defined by the bihomogeneous equation
$u^n v^m E_{P,Q}(x/u,y/v)=0$, $m=\deg Q$.
We identify $\C\P^1$ with $\bC$ denoting $(x:1)$ by $x$ and $(1:0)$ by $\infty$.

Suppose there exists a proper factor $E'(x,y)$ of $E_{P,Q}$. Let $C'$ be
the corresponding subset of $C$. Let $y_1,\dots,y_l\in\C\P^1$ be the poles of $Q$
of multiplicities $q_1,\dots,q_l$ respectively.
The germ of $C$ at $(\infty,y_i)$ has the form $U^n=Y^{q_i}$ in some local
analytic coordinates $(U,Y)$. Indeed, the equation of $C$ in the affine coordinates
$\widehat U=u/x$, $\widehat Y=y/v-y_i$ ($\widehat Y=v/y$ if $y_i=\infty$) has the form
$$
  \widehat U^n f_i(\widehat U)=\widehat Y^{q_i}g_i(\widehat Y),\qquad f_i(0)g_i(0)\ne 0,
$$
thus it has the desired form in the local coordinates
$U=\widehat U f_i^{1/n}$, $Y=\widehat Y g_i^{1/q_i}$ for any choice
of one-valued branches of the roots of $f_i$ and $g_i$ near $0$. 

The binomial $U^n-Y^{q_i}$ factorizes as
$\prod_j^{k_i}(U^{b_i}-\omega^j Y^{a_i})$, where $k_i=\gcd(q_i,n)$, $a_i=q_i/k_i$,
$b_i=n/k_i$, and $\omega$ is a primitive $k_i$-th root of unity. Thus the germ of
$C$ at $(\infty,y_i)$ has $k_i$ local analytic branches, which we denote by
$\gamma_{ij}$, $j=1,\dots,k_i$.
Let $k'_i$ be the number of those that belong to $C'$.

Let $L_i=\C\P^1\times\{y_i\}$. For local intersections, we have
$(\gamma_{ij}\cdot L_i)_{(\infty,y_i)}=b_i$ for each $i,j$. Hence
$$
    k'_i b_i = (C'.L_i)_{(\infty,y_i)}=n' := \deg_x E', \qquad i=1,\dots,l,
$$
and we obtain
\be
  \l{imp}
     \frac{k'_i a_i}{q_i} = \frac{k'_i}{k_i} = \frac{k'_i b_i}{n} = \frac{n'}{n},
     \qquad 1\le i\le l.
\ee
Let $d'/d$ be the reduced form of this fraction, i.e., $d'/d=n'/n$ and $\gcd(d',d)=1$.
Then $d>1$ (because $n'<n$) and $d$ divides all the denominators in \eqref{imp},
which implies that $\gcd(q_1,\dots,q_l,n)>1$.
\qed
 
\medskip

Now we expose more or less the same proof using the language of field extensions.
The notations in both proofs are consistent.
Of course, the proof of \cite[Theorem~1]{o} can be reinterpreted
similarly.
 
\medskip
\noindent {\it Second proof of Theorem \ref{tp}.} 
For a compact Riemann surface $C$, we denote the field of meromorphic functions on $C$
by $\f M(C)$. Given a meromorphic function \hbox{$\theta: C\rightarrow \C\P^1$},
we denote the local multiplicity of $\theta$ at the point $z$ by $e_{\theta}(z)$.

Assume that the curve $E_{P,Q}(x,y)=0$ is reducible, and let $C$ be a desingularization
of some of its components. Then there exist meromorphic functions
$\phi:C\rightarrow \C\P^1$ and $\psi:C\rightarrow \C\P^1$ of degrees $m'<m=\deg Q$
and $n'<n$ such that \be \l{ps} P\circ \phi=Q\circ \psi\ee and
the compositum of the fields $\phi^*\f M(\C\P^1)\subseteq \f M(C)$ and
$\psi^*\f M(\C\P^1)\subseteq \f M(C)$ is the whole field $\f M(C)$.
Furthermore, if  $\theta:C\rightarrow \C\P^1$ is a meromorphic function defined by any
of the sides of equality \eqref{ps}, then by the Abhyankar Lemma
(see e.~g.~\cite{sti}, Theorem 3.9.1) for every point $t_0$ of $C$
the following equality holds: 
\be
   \l{abj} e_{\theta}(t_0)={\rm lcm} \big(e_{P}(\phi(t_0)), e_{Q}(\psi(t_0))\big)
\ee

Let $Q^{-1}\{\infty\}=\{y_1,y_2,\dots,y_l\},$ where $e_Q(y_i)=q_i$, $1\leq i \leq l$, and 
$$
    \psi^{-1}\{y_i\}=\{z_{i1},z_{12},\dots,z_{ik'_i}\},\quad 1\leq i \leq l.
$$
Let us set
$$
     k_i=\gcd(q_i,n), \quad a_i=q_i/k_i, \quad b_i=n/k_i.
$$
Since $P^{-1}\{\infty\}=\infty$, by \eqref{abj} we have 
$$
    e_{\theta}(z_{ij})={\rm lcm}(q_i,n),\quad 1\leq j \leq k'_i,\quad 1\leq i \leq l.
$$  
 Therefore,  
$$
    e_{\psi}(z_{ij}) = \frac{{\rm lcm}(q_i,n)}{q_i} = b_i, \qquad
          1\leq j \leq k'_i,\quad 1\leq i \leq l,
$$ 
implying that, for each $i=1,\dots,l$, we have
$k'_ib_i = \deg\psi = n'$, whence we obtain the equation \eqref{imp} and conclude in
the same way as in the first proof.
\qed


\subsection{A counterexample}\l{sect.ce}
In this subsection we show that the Composition Condition (see Definition~\ref{def.cc})
is not necessary for reducibility of the curve \hbox{$L_P(x,y)=0$}.
Let us set  
\be
    \l{if} T(z) = \frac{z-i}{z+i}.
\ee 
We recall that this M\"obius transformation 
(called the Cayley transform) maps $\bR:=\R\cup\{\infty\}$
to the unit circle $\T$. 

\bl \l{u2}
Let $S$ be a non-constant complex rational function. Then the curve
$\E_{T\circ S}(x,y)=0$ is irreducible if and only if
the curve $E_{S,\overline{S}}(x,y)=0$ is irreducible. 
\el

\pr If $S=S_1/S_2$, where $S_1$ and $S_2$ have no common roots, then
$$
    (T\circ S)(z)=\frac{S_1(z)-iS_2(z)}{S_1(z)+iS_2(z)}.
$$ 
Moreover, it is easy to see that $S_1(z)-iS_2(z)$
and $ S_1(z)+iS_2(z)$ have no common roots. 
Thus,  
\begin{align*}
   \E_{T\circ S}(x,y)&=(S_1(x)-iS_2(x))(\overline{S}_1(y)+i\overline{S}_2(y))
                    -(S_1(x)+iS_2(x))(\overline{S}_1(y)-i\overline{S}_2(y))
\\
    &=2i(S_1(x)\overline{S}_2(y)-S_2(x)\overline{S}_1(y))
    =2iE_{S,\overline{S}}(x,y). \hskip34.8mm{\Box}
\end{align*}

To show that the Composition Condition is not necessary for the reducibility of
$L_P(x,y)$, one can use examples of reducible curves of the form
$E_{S,\overline{S}}(x,y)=0$ found in \cite{cc}. For instance, we can take 
\be
 \l{id}
 \begin{split} 
S(z) &=\frac{1}{11} z^{11} - (a+1)z^9 + 2z^8 + (3a-9)z^7 - 16(a+1)z^6 + (21a+36)z^5\\
 &\;+ (30a-90)z^4 - 63a z^3 + (100a+120)z^2 + (24a-117)z - 18(a+1),
\end{split} 
\ee
where $a$ satisfies $a^2 + a + 3 = 0$.
It is shown in \cite{cc} that for $S(z)$ the following conditions hold.
First, the curve $S(x)-\overline{S}(y)=0$ is reducible. Second,  
\be
    \l{cont}
              \overline{S}(z)\neq S(cz+b),
\ee 
for any  $c\in \C^*$, $b\in \C$ (the last condition makes this example
non-trivial since every  curve of the form 
$S(x)-S(cy+b)=0$ obviously has a factor $x-cy-b=0$).  

Let us consider a rational function $P=T\circ S$, where $T$ and $S$ are defined by
\eqref{if} and \eqref{id}.   
By Lemma~\ref{u2}, the  curve $\E_P(x,y)=0$ is reducible, implying by Lemma~\ref{u1}
that the  curve $L_P(x,y)=0$ is also reducible. 
Since the degree of $P$ is a prime number, if
the Composition Condition were satisfied by $P$,
we would have
\be
   \l{iwi}
              P=T\circ S=B\circ \mu
\ee
for some quotient of finite Blaschke products $B$ and M\"obius transformation $\mu$.
Let now $\nu$ be a M\"obius transformation that maps
$\bR$ to $\mu^{-1}(\T)$. Then $(\mu\circ \nu)(\bR)=\T$ and it follows from \eqref{iwi}
combined with the equalities $B(\T)=\T$ and $T(\bR)=\mathbb T$ that 
\be
        (S\circ \nu)(\bR)=(S\circ\mu^{-1})(\T)=(T^{-1}\circ B)(\T)=T^{-1}(\T)=\bR.
\ee 
Therefore the rational function $\overline{(S\circ\nu)}-(S\circ\nu)$ identically vanishes on $\R$,
and hence  on $\C$. Thus,
$$
    S\circ \nu=\overline{S\circ \nu}=\overline{S}\circ \overline{\nu},
$$
and denoting the M\"obius transformation $\nu\circ \overline{\nu}^{-1}$ by $\delta$,
we rewrite this equality as
\be
   \l{eq.last}
              \overline{S}=S\circ \delta.
\ee
Since $S$ is a polynomial, \eqref{eq.last} implies that $\delta$ also  
is a polynomial. Moreover, \hbox{$\deg\delta=1$} (because $\delta$ is a M\"obius
transformation), thus \eqref{eq.last} contradicts \eqref{cont}.
The obtained contradiction shows that $P$ does not satisfy the Composition Condition.


\begin{thebibliography}{10}

\bibitem{AR} N.~Ailon and Z.~Rudnick,  \textit{Torsion points on curves and common divisors of $a^k-1$ and $b^k-1$},  Acta Arith.,  {\bf 113} (2004),   31--38. 

\bibitem{b} D.\,N.~Bernshtein, \textit{The number of roots of a system of equations},
Funct. Anal. Appl., {\bf 9}:3 (1975), 183--185.

\bibitem{c} M.\,L.~Cartwright, \textit{On the level curves of integral and meromorphic functions}, Proc. London Math.
Soc., {\bf  43} (1937), 468--474.

\bibitem{cc} P.~Cassou-Nogu\`es, J.-M.~Couveignes,
\textit{Factorisations explicites de $g(y)-h(z)$},
Acta Arith. {\bf 87} (1999), 291--317. 

\bibitem{CMZ}
P.~Corvaja, D.~Masser and U.~Zannier,
\textit{Sharpening `Manin--Mumford' for certain algebraic groups of dimension 2}, 
Enseign. Math., {\bf 59} (2013), 1--45.

\bibitem {f4} M.~Fried, \textit{
The field of definition of function fields and a problem in the reducibility of polynomials in two variables}, Ill. J. Math. {\bf 17} (1973), 128--146.

\bibitem {fc} M.~Fried, {\it  Variables separated equations: strikingly different roles
for the branch cycle lemma and the finite simple group classification},
Sci. China Math. {\bf 55} (2012), no. 1, 1--72.

\bibitem {ful} W.~Fulton, {\it Intersection theory},
Ergeb. Math. Grenzgeb. (3), Bd. 2, Springer-Verlag, Berlin, 1984.

\bibitem {gh} \'E.~Ghys, \textit{A singular mathematical promenade},
ENS Éditions, Lyon, 2017.

\bibitem{nef2} J.~K\"onig, D.~Neftin,  \textit{Reducible Fibers of Polynomial Maps},
arXiv:2001.03630.  

\bibitem{ok} Y.~Okuyama, M.~Stawiska, \textit{A generalization of the converse of
Brolin's theorem}, Ann. Polon. Math. {\bf 123} (2019), no.~1, 467--472.

\bibitem{o} S.\,Yu.~Orevkov, {\it Irreducibility of lemniscates},
Russian Math. Surveys {\bf 73} (2018), no. 3, 543--545.

\bibitem{ps} F.~Pakovich, I.~Shparlinski, {\it Level curves of rational functions and
unimodular points on rational curves,} Proc. Amer. Math. Soc. {\bf 148} (2020), no. 5,
1829--1833.

\bibitem{shaf} I.\,R.~Shafarevich, {\it Basic algebraic geometry. 1. Varieties in projective space.} Third edition. Translated from the 2007 third Russian edition. Springer, Heidelberg, 2013.

\bibitem{stsu} K.~Stephenson and  C.~Sundberg, \textit{Level curves of inner functions}, Proc. London Math. Soc., {\bf  51} (1985),   77--94.

\bibitem{ste} K.~Stephenson, \textit{Analytic functions sharing level curves and tracts}, Ann. of Math., {\bf 123} (1986),   107--144.

\bibitem {sti} H.~Stichtenoth, \textit{Algebraic Function Fields and Codes,} second ed.,
Graduate Textbooks in Mathematics 254, Springer-Verlag, Berlin, 2009.

\end{thebibliography}
\end{document}